# Vulnerability Analysis for Complex Networks Using Aggressive Abstraction

Richard Colbaugh     Kristin Glass

*Abstract*—Large, complex networks are ubiquitous in nature and society, and there is great interest in developing rigorous, scalable methods for identifying and characterizing their vulnerabilities. This paper presents an approach for analyzing the dynamics of complex networks in which the network of interest is first abstracted to a much simpler, but mathematically equivalent, representation, the required analysis is performed on the abstraction, and analytic conclusions are then mapped back to the original network and interpreted there. We begin by identifying a broad and important class of complex networks which admit *vulnerability-preserving, finite state abstractions*, and develop efficient algorithms for computing these abstractions. We then propose a vulnerability analysis methodology which combines these finite state abstractions with formal analytics from theoretical computer science to yield a comprehensive vulnerability analysis process for networks of real-world scale and complexity. The potential of the proposed approach is illustrated with a case study involving a realistic electric power grid model and also with brief discussions of biological and social network examples.

## I. INTRODUCTION

It is widely recognized that technological, biological, and social networks, while impressively robust in most circumstances, can fail catastrophically in response to focused attacks. Indeed, this combination of robustness and fragility appears to be an inherent property of complex, evolving networks ranging from the Internet and electric power grids to gene regulatory networks and financial markets [e.g., 1-5]. As a consequence, there is significant interest in developing methods for reliably detecting and characterizing the vulnerabilities of these networks [e.g., 6,7].

The challenges of vulnerability analysis are particularly daunting in the case of complex networks. Most such networks are large-scale "systems of systems", so that analysis methods must be computationally efficient. Additionally, because these networks perform reliably almost all of the time, standard techniques for finding vulnerabilities (e.g., computer simulations, "red teaming") can be ineffective and, in any case, are not guaranteed to identify all vulnerabilities. These observations suggest that, in order to be practically useful, any method for analyzing vulnerabilities of complex networks should be *scalable*, to enable analysis of networks of real-world complexity, and *rigorous*, so that for instance it is guaranteed to find all vulnerabilities of a given class.

This paper presents a new approach to vulnerability analysis which possesses these properties. The proposed methodology is based upon *aggressive abstraction* – dramatically simplifying, property preserving abstraction of the network of interest [4]. Once an aggressive abstraction is derived, all required analysis is performed using the abstraction. Analytic conclusions are then mapped back to the original network and interpreted there; this mapping is possible because of the property preserving nature of the abstraction procedure.

Our focus is on dynamical systems with uncountable state spaces, as many complex networks are of this type. We begin by identifying a large and important class of dynamical networks which admit *vulnerability-preserving, finite state abstractions*, and develop efficient algorithms for recognizing such networks and for computing their abstraction. We then offer a methodology which combines these finite state models with formal analytics from theoretical computer science [8] to provide a comprehensive vulnerability analysis process for networks of real-world scale and complexity. The potential of the proposed approach to complex networks analysis is illustrated through a case study involving vulnerability analysis of a realistic electric power grid and also via brief discussions of biological and social network examples.

## II. PRELIMINARIES

This section introduces the class of network models to be considered in the paper and briefly summarizes some technical background that will be useful in our development.

The evolution to ensure robust performance in complex networks typically leads to systems that possess a "hybrid" structure, exhibiting both continuous and discrete dynamics [4]. More precisely, these networks often evolve to become *hybrid dynamical systems* – feedback interconnections of switching systems, which have discrete state sets, with systems whose dynamics evolve on continuous state spaces [9].

More quantitatively, consider the following definitions for hybrid dynamical system (HDS) models:

**Definition 2.1:** A *continuous-time HDS* is a control system

$$\Sigma_{HDSct} \quad \begin{aligned} q+ &= h(q,k), \\ dx/dt &= f_q(x,u), \\ k &= p(x), \end{aligned}$$

where $q \in Q$ (with $|Q|$ finite) and $x \in X \subseteq \Re^n$ (with X bounded) are the states of the discrete and continuous systems that

The research described in this paper was supported by the U.S. Department of Defense, the U.S. Department of Energy, the U.S. Department of Homeland Security, and the Laboratory Directed Research and Development program at Sandia National Laboratories.

R. Colbaugh is with Sandia National Laboratories, Albuquerque, NM 87111 USA and New Mexico Institute of Mining and Technology, Socorro, NM 87801 USA (phone: 505-603-1248; e-mail: colbaugh@nmt.edu).
K. Glass is with New Mexico Institute of Mining and Technology, Socorro, NM 87801 USA (e-mail: kglass@icasa.nmt.edu).

make up the HDS, u∈$\Re^m$ is the control input, h defines the discrete system dynamics, $\{f_q\}$ is a family of vector fields characterizing the continuous system dynamics, and p defines a partition of state space X into subsets with labels k∈{1, …, K}.

**Definition 2.2:** A *discrete-time HDS* is a control system

$$\Sigma_{HDSdt} \quad \begin{aligned} q+ &= h(q,k), \\ x+ &= f_q(x,u), \\ k &= p(x). \end{aligned}$$

We sometimes refer to an HDS using the symbol $\Sigma_{HDS}$ if the nature of the continuous system (continuous- or discrete-time) is either unimportant or clear from the context.

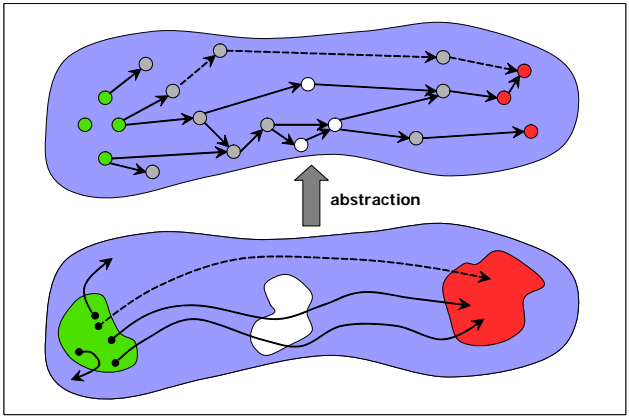

Figure 1. Finite state abstraction. Cartoon illustrates that abstraction preserves network dynamics: trajectories of the infinite state system (curves in blue region at bottom) are mapped to equivalent finite state trajectories (sequences of state transitions at top).

The concept of finite state abstraction for an infinite state system is illustrated in Figure 1. Consider a complex network with states that evolve on a continuous space and an analysis question of interest. Such a situation is depicted at the bottom of Figure 1, where the continuous dynamics are shown as curves on a continuous state space (blue region), and the analysis question involves deciding whether states in the green region can evolve to the red region. Reachability questions of this sort are difficult to answer for generic complex networks. However, if it is possible to construct a finite state abstraction of the network which possesses equivalent dynamics, then the analysis task becomes much easier. To see this, observe from Figure 1 that a finite state abstraction of the original dynamics takes the form of a graph, where the states are graph vertices (nodes within the blue region at top) and feasible state transitions define the graph's directed edges. Reachability analysis is straightforward with a graph, and if the complex network and its abstraction have equivalent reachability properties then the much simpler graph analysis also characterizes the reachability of the original system.

Reachability assessment, while valuable, is typically not sufficient to answer real-world vulnerability analysis questions. For instance, suppose that the red region in Figure 1 is the set of failure states. It may be of interest to determine if all system trajectories which reach the red region first pass through the white "alerting" region, so there is warning of impending failure, or whether all trajectories which reach the red region subsequently return to the blue "normal" region, and thereby recover from failure. Addressing these more intricate questions requires that the analysis be conducted using a language which allows a nuanced description of, and reasoning about, network dynamics. We show in [4] that *linear temporal logic* (LTL) provides such a language, enabling quantitative specification of all vulnerability problems we've encountered in complex networks analysis. LTL extends propositional logic by including temporal operators, thereby allowing dynamical phenomena to be analyzed, and is similar to natural language and thus is easy to use [10].

As we wish to use LTL to analyze the dynamics of complex networks and we model these networks as HDS, we tailor our definition of LTL to be compatible with this setting:

**Definition 2.3:** The *syntax* of LTL consists of

- *atomic propositions* (q,k), where q ∈ Q is an HDS discrete state and k ∈ K is a label for a subset in the continuous system state space partition;

- *formulas* composed from atomic propositions using a grammar of Boolean (φ ∨ θ, ¬φ) and temporal (φ**U**θ, ○φ) operators.

The *semantics* of LTL follows from interpreting formulas on trajectories of HDS, that is, on sequences of (q, k) pairs: (**q**, **k**) = $(q_0, k_0), (q_1, k_1), \ldots, (q_T, k_T)$.

The Boolean operators ∨ and ¬ are disjunction and negation, as usual. The temporal operators **U** and ○ are read "until" and "next", respectively, with φ**U**θ specifying that φ must hold until θ holds and ○φ signifying that φ will be true at the next time instant (see [10] for a more thorough description).

Abstractions which preserve LTL also preserve vulnerabilities [4]. Thus we seek an abstraction procedure which preserves LTL: given a system representation $\Sigma_1$, the procedure should generate a system abstraction $\Sigma_2$ which is such that $\{\Sigma_1 \models \varphi\} \Leftrightarrow \{\Sigma_2 \models \varphi\}$ for all LTL formulas φ (where $\models$ denotes formula satisfaction). Bisimulation is a powerful method for abstracting *finite* state systems to yield simpler finite state systems that are equivalent from the perspective of LTL [10]. However, the problem of constructing finite state bisimulations for continuous state systems is largely unexplored (but see the seminal work [11,12]). Indeed, one of the contributions of this paper is to develop a theoretically sound, practically implementable approach to obtaining finite state bisimulations for complex network models.

Bisimulation is typically defined for transition systems, so we first introduce this notion (see [10] for details):

**Definition 2.4:** A *transition system* is a four-tuple T = (S, →, Y, h) with state set S, transition relation →⊆ S × S, output set Y, and output map h: S → Y. T is *finite* if |S| is finite.

Transition relation → defines admissible state transitions, so (q, q')∈→, denoted q → q', if T can transition from q to q'.

Bisimilar transition systems share a common output set and have dynamics which are equivalent from the perspective of these outputs:

**Definition 2.5:** Transition systems $T_S = (S, \to_S, Y, h_S)$ and $T_P = (P, \to_P, Y, h_P)$ are *bisimilar* via relation $R \subseteq S \times P$ iff:

- s ~ p ⇒ $h_S(s) = h_P(p)$ (R respects observations);
- s ~ p, s $\to_S$ s' ⇒ ∃ p' ~ s' such that p $\to_P$ p' ($T_P$ simulates $T_S$, denoted $T_S \angle T_P$);
- p ~ s, p $\to_P$ p' ⇒ ∃ s' ~ p' such that s $\to_S$ s' ($T_P \angle T_S$).

A standard result from theoretical computer science [e.g., 10] shows that bisimulation preserves LTL:

**Proposition 2.1:** If $T_1$ and $T_2$ are bisimilar transition systems and φ is an LTL formula then $\{T_1 \models \varphi\} \Leftrightarrow \{T_2 \models \varphi\}$.

The following alternative definition for bisimulation is easily shown to be equivalent to the one presented in Definition 2.5 and is useful in the subsequent development:

**Definition 2.6:** A finite partition $\Phi: S \to P$ of the state space S of transition system $T = (S, \to, Y, h)$ naturally induces a *quotient transition system* $T/\sim = (P, \to_\sim, Y, h_\sim)$ of T provided that

- $\Phi(s) = \Phi(s')$ (denoted s ~ s') ⇒ h(s) = h(s');
- $h_\sim(p) = h(s)$ if p = $\Phi(s)$;
- $\to_\sim$ is defined so that $\Phi(s) \to_\sim \Phi(s')$ iff s → s'.

Transition system T and its quotient T/~ are bisimilar if an additional condition holds:

**Proposition 2.2:** Suppose T/~ is defined as in Definition 2.6 and, in addition, $\Phi(s) \to_\sim \Phi(s') \Rightarrow \forall \, s'' \sim s \, \exists \, s''' \sim s'$ such that s'' → s'''. Then T and T/~ are bisimilar.

Finally, we introduce a class of continuous state (control) systems which is important in applications.

**Definition 2.7:** The continuous-time system dx/dt=f(x,u), with f: $\Re^n \times \Re^m \to \Re^n$, is *differentially flat* if there exists (flat) outputs $z \in \Re^m$ such that z = H(x), x = $F_1(z, dz/dt, …, d^r z/dt^r)$, and u = $F_2(z, dz/dt, …, d^r z/dt^r)$ for some integer r and maps H, $F_1$, $F_2$.

**Definition 2.8:** The discrete-time system x+ = f(x,u) is *difference flat* (with memory k) if there exists (flat) outputs $z \in \Re^m$ such that z = H(x), x(t) = $F_1(z(t), z(t+1), …, z(t+k−1))$, and u(t) = $F_2(z(t), z(t+1), …, z(t+k−1))$ for some maps H, $F_1$, $F_2$.

Background on flat systems may be found in [13]. Many real-world control systems are flat, including all controllable linear systems as well as all feedback linearizable systems. Perhaps more importantly, the complex, evolving networks underlying so much of advanced technology, biology, and social processes frequently possess flat subsystems.

## III. FINITE STATE ABSTRACTION

In this section we demonstrate that hybrid systems with (differentially or difference) flat continuous systems admit finite state bisimulations and present algorithms for constructing the bisimilar abstractions.

Consider an HDS of the form given in Definition 2.1 or 2.2. The following provides a transition system representation for the continuous system dynamics of HDS:

**Definition 3.1:** The *transition system model* $T_{HDSc}$ for the continuous system portion of $\Sigma_{HDS}$ is the collection $T_{HDSc} = \{T_q^k\}$, with one transition system $T_q^k = (X_q^k, \to_q^k, Y_q^k, h_q^k)$ specified for each (q, k) pair. Each $T_q^k$ has bounded state space $X_q^k$, finite output set $Y_q^k$, an output map $h_q^k: X_q^k \to Y_q^k$ that defines a finite partition of $X_q^k$ with labels $y \in Y_q^k$, and transition relation $\to_q^k$ reflecting the discrete- or continuous-time dynamics:

- for discrete-time continuous systems, x $\to_q^k$ x' iff ∃u such that x' = $f_q(x,u)$ on subset k;
- for continuous-time continuous systems, x $\to_q^k$ x' iff there is a trajectory x: [0, T] → $X_q^k$ of dx/dt = $f_q(x,u)$, a time t'∈ (0,T), and adjacent partitions of $X_q^k$ labeled y, y'∈$Y_q^k$ such that x(0)=x, x(T)=x', x([0,t'))⊆y, and x((t',T])⊆y'.

We make the standard assumption that k: X → K partitions the HDS continuous system state space X into polytopes and that all HDS discrete system transitions are triggered by k transitions [9] (see Definitions 2.1 and 2.2).

Definition 3.1 allows $\Sigma_{HDS}$ to be modeled as a feedback interconnection of two transition systems, one with continuous state space and one with finite state set:

**Definition 3.2:** The *transition system* $T_{HDS}$ associated with the HDS given in Definition 2.1 or 2.2 is a feedback interconnection of 1.) the continuous system transition system $T_{HDSc} = \{T_q^k\}$ given in Definition 3.1 and 2.) the transition system associated with the HDS discrete system, given by $T_{HDSd} = (Q, \to_d, Q, id)$, where id is the identity map and q $\to_d$ q' iff ∃k such that q' = h(q, k). Thus $T_{HDS} = (Q \times X, \to_{HDS}, Q \times Y, h_{HDS})$, where $Q \times X = \cup_q (\cup_k \{q\} \times X_q^k)$, $Q \times Y = \cup_q (\cup_k \{q\} \times Y_q^k)$, and the definitions for $\to_{HDS}$ and $h_{HDS}$ follow immediately from the transition relation and output map definitions specified for $T_{HDSc}$ and $T_{HDSd}$.

Because the transition system $T_{HDSd}$ corresponding to the HDS discrete system is already a finite state system, the main challenge in abstracting HDS to finite state systems is associated with finding finite state bisimulations for the continuous systems $T_{HDSc} = \{T_q^k\}$. This is made explicit in the following

**Theorem 1**: If each transition system $T_q^k$ associated with $T_{HDS}$ is bisimilar to its finite quotient transition system $T_q^k/\sim = (Y_q^k, \to_\sim, Y_q^k, id)$ and the state space quotient partitions defined by the $h_q^k$ satisfy a mild compatibility condition then $T_{HDS}$ admits a finite bisimulation.

**Proof:** The proof is straightforward and is given in [4]. ∎

Theorem 1 shows that the key step in obtaining a finite state bisimulation for HDS $T_{HDS}$, and thus for $\Sigma_{HDS}$, is constructing bisimulations for the continuous state transition systems $T_q^k$. We therefore focus on this latter problem for the remainder of the section. Our first main result along these lines is for difference flat continuous systems:

**Theorem 2:** Given any finite partition $\pi: Z \to Y$ of the flat output space Z of a difference flat system, the associated transition system $T_F = (X, \to, Y, \pi \circ H)$ admits a bisimilar quotient $T_F/\sim$.

**Proof:** Consider the equivalence relation R that identifies state pairs $(x, x')$ which generate identical sets of k-length output symbol sequences $y = y_0 y_1 \ldots y_{k-1}$, and the quotient system $T_F/\sim$ induced by R. R defines a finite partition of X (both |Y| and k are finite), and $x \sim x' \Rightarrow \pi \circ H(x) = \pi \circ H(x')$ so that R respects observations. $T_F \angle T_F/\sim$ follows immediately from the definition of quotient systems. To see that $T_F/\sim \angle T_F$, note that flatness ensures *any* symbol string $y = y_k y_{k+1} \ldots$ is realizable by transition system $T_F$; thus $x \sim x'$ at time t implies that x and x' can transition to equivalent states at time t + 1. Therefore, from Definition 2.5, $T_F$ and $T_F/\sim$ are bisimilar. ∎

**Remark 3.1:** Efficient algorithms exist for checking if a given system is difference flat, so Theorem 2 provides a practically implementable means of identifying discrete-time continuous state systems which admit finite bisimulation [4].

**Remark 3.2:** The flat output trajectory completely defines the evolution of a difference flat system. As a consequence, because *any* finite partition of flat output space induces a finite bisimilar quotient for the flat system, this partition can be refined to yield any desired level of detail in the abstraction.

An analogous result holds for differentially flat HDS continuous systems. Our development of this result requires the following lemmas.

**Lemma 3.1:** A control system is differentially flat iff it is dynamic feedback linearizable.
**Proof:** The proof is given in [14]. ∎

**Lemma 3.2:** Control system $\Sigma$ admits a finite bisimulation iff any representation of $\Sigma$ obtained through coordinate transformation and/or invertible feedback also admits a finite bisimulation.
**Proof:** The proof is straightforward. ∎

Lemmas 3.1 and 3.2 suggest the following procedure for constructing finite bisimulations for differentially flat systems: 1.) transform the flat system into a linear control system via feedback linearization, 2.) compute a finite bisimulation for the linear system, and 3.) map the bisimilar model back to the original system representation. As a result, we focus on building finite bisimulations for linear control systems.

In particular, the control system of interest is one "chain" of a Brunovsky normal form (BNF) system $\Sigma_{BNF}$ [4]:

$$dx_1/dt = x_2,$$
$$dx_2/dt = x_3,$$
$$\ldots$$
$$dx_n/dt = u.$$

Concentrating on this system entails no loss of generality, as any controllable linear system can be modeled as a collection of these single chain systems, one for each input, and the decoupled nature of the chains ensures we can abstract each one independently and then "patch" the abstractions together to obtain an abstraction for the full system.

Consider the following partition of the (assumed bounded) state space $X \subseteq \Re^n$ of $\Sigma_{BNF}$:

**Definition 3.3:** *Partition* $\pi_\varepsilon$ is the map $\pi_\varepsilon: X \to Y$ that partitions X into subsets $y_{is} = \{x \in X \mid x_1 \in [i\varepsilon, (i+1)\varepsilon), \text{sign}(x_2) = s_1, \ldots, \text{sign}(x_n) = s_{n-1}\}$, where i is an integer and s is an (n–1)-vector of "signs" specifying a particular orthant of X.

Note that $\pi_\varepsilon$ partitions X into "slices" orthogonal to the $x_1$-axis.

We are now in a position to state

**Theorem 3:** The transition system $T_{BNF} = (X, \to, Y, \pi_\varepsilon)$ associated with system $\Sigma_{BNF}$ and partition $\pi_\varepsilon$ admits a finite bisimilar quotient $T_{BNF}/\sim = (Y, \to_\sim, Y, id)$.

**Proof:** $T_{BNF}/\sim$ is finite because |Y| is finite. Assume $\to_\sim$ is constructed so that $\pi_\varepsilon(x) \to_\sim \pi_\varepsilon(x') \Leftrightarrow x \to x'$, with the latter specified as in Definition 3.1. Then all the conditions of Definition 2.6 are satisfied, and from Proposition 2.2 we need only show $\pi_\varepsilon(x) \to_\sim \pi_\varepsilon(x') \Rightarrow \forall x'' \sim x \exists x''' \sim x'$ such that $x'' \to x'''$. This amounts to demonstrating that if x can be driven through face F of slice $\pi_\varepsilon(x)$ then any $x''$ in that slice can be driven through F as well, which can be shown by checking this property for the system $x_1^{(n)} = u$ on each orthant of X. ∎

**Remark 3.3:** The result given in Theorem 3 is most useful in situations where the control input u can be chosen large relative to the system "drift". Abstraction methods for applications in which control authority is limited are given in [4].

Next we turn to the task of *computing* finite bisimulations for HDS. We focus on constructing bisimulations for HDS continuous systems, as HDS discrete systems already possess finite state representations, and in particular on abstracting differentially flat continuous systems; the derivation of algorithms for difference flat continuous systems is analogous but simpler and is therefore omitted (see [4]).

Consider, without loss of generality (see Lemma 3.1), the problem of computing a finite state abstraction for continuous-time linear control system $\Sigma_{lc}: dx/dt = Ax + Bu$ (where A,B are matrices). The transition system associated with $\Sigma_{lc}$ is $T_{lc} = (X, \to_{lc}, Y, h)$, with $h: X \to Y$ any finite, hypercubic partition of X and $\to_{lc}$ specified as in Definition 3.1. The finite state abstraction of interest is quotient system $T_{lc}/\sim = (Y, \to_{lc\sim}, Y, id)$. Observe that in order to obtain $T_{lc}/\sim$ it is only necessary to determine the set of admissible transition relations $\to_{lc\sim}$.

As most applications of interest involve large-scale systems, it is desirable to develop efficient algorithms for computing $\to_{lc\sim}$. We now introduce such a procedure. The algorithm decides whether a transition $y \to_{lc\sim} y'$ between two adjacent cells of the lattice y, y' is allowed, and is repeated for all candidate transitions of interest. We begin by summarizing a simple algorithm, based on computational linear system results given in [15], for deciding whether $y \to_{lc\sim} y'$ is admissible. Let k be the number of the coordinate axis orthogonal to the common face between y and y', V be the set of vertices shared by y and y', and $a_k^T$ represent row k of A. Define $\Pi^k(w)$ to be the projection of vector w onto axis k, and suppose $y < y'$. Then $y \to_{lc\sim} y'$ iff $\Pi^k(Av_i + Bu) > 0$ for some $v_i \in V$ and $u \in U$. An algorithm which "operationalizes" this observation is

**Algorithm 3.1:**

If $y < y'$:

- If any element of row k of B is nonzero, $y \to_{lc\sim} y'$ is true. STOP.
- Repeat until $y \to_{lc\sim} y'$ is determined to be true or all vertices have been checked:
  o Select a vertex $v_i \in V$.
  o Compute the inner product $p = a_k^T v_i$.
  o If $p > 0$ then $y \to_{lc\sim} y'$ is true. STOP.
- If $y \to_{lc\sim} y'$ has not been found to be true it is false.

If $y > y'$: Algorithm is the same except that the comparison $p > 0$ is replaced by $p < 0$.

A difficulty with Algorithm 3.1 is that the number of vertices shared by two adjacent cells is $2^{n-1}$, so that checking them becomes unmanageable even for moderately-sized systems. Interestingly, the algorithm can be modified so that feasibility of a transition can be tested by considering only a *single* well-chosen vertex, independent of the size of the model [16]. The new algorithm is therefore extremely efficient and can be applied to very large systems. Let $v_0$ be the lowest vertex (in a component-wise sense) shared by y, y' and let $a_k^+$ ($a_k^-$) be the sum of positive (negative) elements of row k of A, excluding the diagonal. We can now state

**Algorithm 3.2** [16]**:**

If $y < y'$:

- If any element of row k of B is nonzero, $y \to_{lc\sim} y'$ is true. STOP.
- Compute the inner product $p = a_k^T v_0$.
- If $p + a_k^+ > 0$ then $y \to_{lc\sim} y'$ is true. STOP.
- Otherwise $y \to_{lc\sim} y'$ is false.

If $y > y'$: Algorithm is the same except that the comparison $p + a_k^+ > 0$ is replaced by $p + a_k^- < 0$.

A Matlab program which implements Algorithm 3.2 is presented in [4]. This program has been applied to systems with n = 10 000 state variables using desktop computers.

## IV. VULNERABILITY ANALYSIS

This section considers the *vulnerability assessment* problem: given a complex network and a class of failures of interest, does there exist an attack which causes the system to experience such failure? Other important vulnerability analysis tasks, including vulnerability exploitation and mitigation, are investigated in [4]. The proposed approach to vulnerability assessment leverages the finite state abstraction results derived in the preceding section. The basic idea is straightforward: given an HDS model for a network of interest and a class of failures of concern: 1.) construct a finite bisimulation for the HDS network model, 2.) conduct the vulnerability analysis on the system abstraction, and 3.) map the analysis results back to the original system model.

Observe that the proposed approach possesses desirable characteristics. For instance, the analytic process is scalable, because both the abstraction methodology and the tools available for detecting vulnerabilities in finite state systems [e.g., 8] are computationally efficient. Additionally, the analysis is rigorous. Because HDS vulnerabilities are expressible as LTL formulas, and bisimulation preserves LTL, the original complex network and its abstraction have identical vulnerabilities. Formal analysis tools such as model checking [8] can be structured to identify all vulnerabilities of the finite state abstraction, and bisimilarity then implies that the approach is guaranteed to find all vulnerabilities of the original network as well.

We now quantify the proposed approach to vulnerability assessment. It is supposed that the complex network of interest can be modeled as an HDS, $\Sigma_{HDS}$, and that the network's desired or "normal" behavior can be characterized with an LTL formula $\varphi$; generalizing the situation to a set of LTL formulas $\{\varphi_i\}$ is straightforward. Consider the following

**Definition 4.1:** Given an HDS $\Sigma_{HDS}$ and an LTL encoding $\varphi$ of the desired network behavior, the *vulnerability assessment problem* involves determining whether $\Sigma_{HDS}$ can be made to violate $\varphi$.

The proposed vulnerability assessment method employs *bounded model checking* (BMC), a powerful technique for deciding whether a given finite state transition system satisfies a particular LTL specification over a finite, user-specified time horizon [8]. Briefly, BMC checks whether a finite transition system T satisfies an LTL specification $\varphi$ on a time interval [0, k], denoted $T \models_k \varphi$, in two steps: 1.) translate $T \models_k \varphi$ to a proposition $[T, \varphi]_k$ which is satisfied by, and only by, transition system trajectories that *violate* $\varphi$ (this is always possible), and 2.) check if $[T, \varphi]_k$ is satisfiable using a modern SAT solver [8]. Note that because modern SAT solvers are extremely powerful, this approach to model checking can be implemented with problems of real-world scale.

We are now in a position to state our vulnerability assessment algorithm. Let $T_{HDS}$ denote the transition system associated with $\Sigma_{HDS}$, and consider the vulnerability assessment problem given in Definition 4.1. We have

**Algorithm 4.1: Vulnerability assessment**

1. Construct a finite bisimilar abstraction T for $T_{HDS}$ using the results of Section 3.
2. Check satisfiability of $[T, \varphi]_k$ using BMC:
   - if $[T, \varphi]_k$ is not satisfiable then T is not vulnerable and thus $\Sigma_{HDS}$ is not vulnerable (on time horizon k);
   - if $[T, \varphi]_k$ is satisfiable then T, and therefore $\Sigma_{HDS}$, is vulnerable, and the SAT solver "witness" is an exploitation of the vulnerability.

To illustrate the utility of the proposed approach to vulnerability assessment we apply the analytic method to an important complex network: an electric power (EP) grid. EP grids are naturally represented as HDS, with the continuous system modeling the generator and load dynamics as well as power flow constraints and the discrete system capturing protection logic switching and other "supervisory" behavior:

**Definition 4.2:** The *HDS power grid model* $\Sigma_{EP}$ takes the form

$$\Sigma_{EP} \quad \begin{aligned} q+ &= h(q, k, v), \\ dx/dt &= f_q(x, y, u), \\ 0 &= g_q(x, y), \\ k &= p(x, y), \end{aligned}$$

where q and x are the discrete and continuous system states, v and u denote exogenous inputs, y is the vector of "algebraic variables", and all other terms are analogous to those introduced in Definition 2.1.

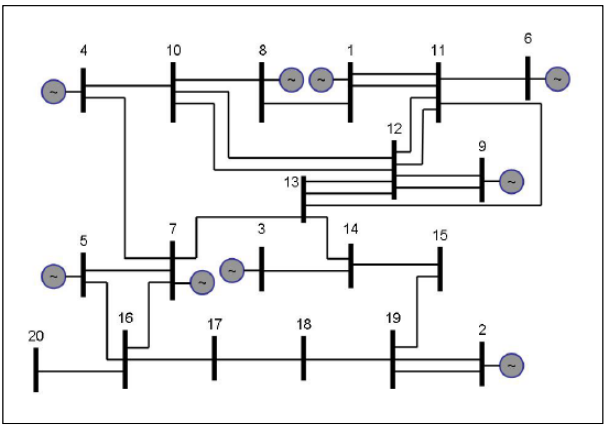

Figure 2. One-line diagram for the 20-bus EP grid model used in the vulnerability assessment case study.

The continuous system portion of grid model $\Sigma_{EP}$ is feedback linearizable [17], which implies the continuous system is differentially flat and consequently that $\Sigma_{EP}$ admits a finite abstraction. Additionally, it can be shown that grid vulnerabilities are expressible as LTL formulas composed of atomic propositions which depend only on q and k [18]. Thus Algorithm 4.1 is directly applicable to power grids.

We now summarize the results of a vulnerability assessment for the 20bus grid shown in Figure 2. This grid provides a simple but useful representation of a real EP system for which (proprietary) data are available to us [4]. The grid can be modeled as an HDS $\Sigma_{EP}$ of the form given in Definition 4.2. The report [4] gives a Matlab encoding of the specific HDS model used in this study. Because the model $\Sigma_{EP}$ corresponds to a real world grid, the behaviors of the model and the actual grid can be compared. For example, the real grid recently experienced a large cascading voltage collapse, and data was collected for this event. We simulated this cascading outage (see Figure 3, top plot) and found close agreement between the behavior of the actual grid and the model $\Sigma_{EP}$. Observe that this result is encouraging given the well-known difficulties associated with reproducing such cascading dynamics with computer models (see, e.g., [17,18]).

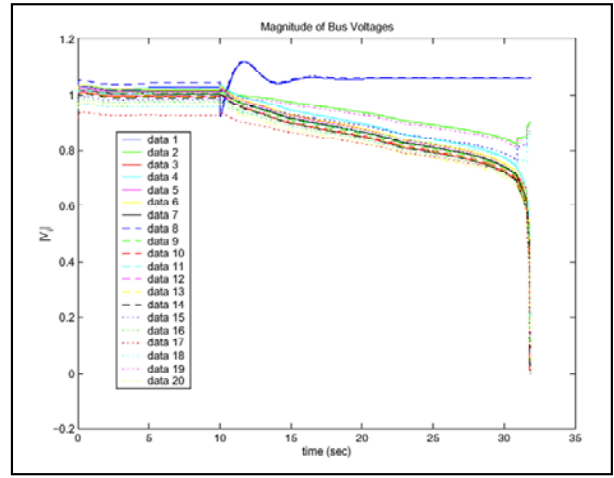

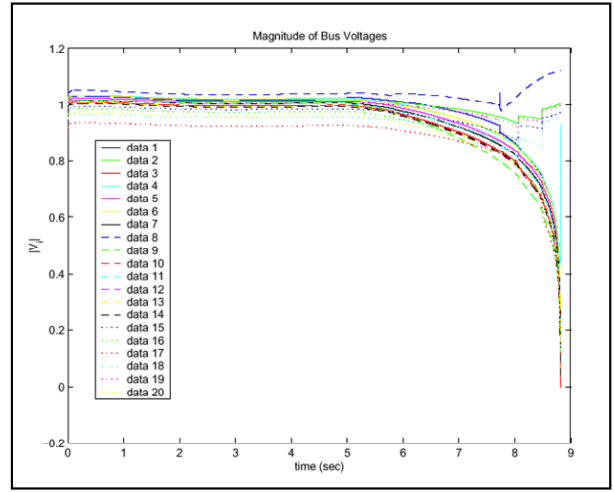

Figure 3. Sample simulation results for 20-bus EP grid model. Plot at top is from the model validation study and shows the evolution of voltages at all 20 buses; these voltage time series are in good agreement with those observed in the corresponding cascading voltage collapse for the actual grid. Plot at bottom depicts voltage time series which result from applying the vulnerability exploitation procedure designed using the proposed finite state abstraction methodology.

Vulnerability assessment was performed using Algorithm 4.1. It was assumed that the grid's attacker wishes to drive the voltage at bus 11 to unacceptably low levels, so that the loads at this bus would not be served, and that the attacker has only limited grid access. In particular, we consider here a scenario in which the attacker can gain assess to the generator at bus 2 via cyber means [18]. Note that this class of vulnerabilities is interesting because the access point – the generator at bus 2 – is geographically remote from the target of the attack – the loads at bus 11.

The first step in the vulnerability assessment procedure specified in Algorithm 4.1 involves constructing a finite state bisimulation T for $\Sigma_{EP}$; this abstraction is computed using Algorithm 3.2. The second step in Algorithm 4.1 is to apply BMC to T to determine if it is possible to realize the attack objective, i.e., low voltage at bus 11, through admissible manipulation of the generator at bus 2. We employed NuSMV, an open source software tool for formal verification of finite state systems, for this analysis [19]. This vulnerability assessment reveals that it is possible for the attacker to realize the given objective via the assumed grid access, and gives a finite state "trace" of one means of exploiting the vulnerability. Using this trace, we synthesize an exploitation attack which is directly implementable with the HDS model $\Sigma_{EP}$. Sample simulation results are shown in Figure 3 (bottom plot). It can be seen from the bus voltage time series in Figure 3 that the attacker's goals can indeed be realized, in this case by initiating a cascading voltage collapse which takes down bus 11 as well as most of the rest of the grid.

## V. Discussion

This paper presents an approach for analyzing complex networks in which the network of interest is first abstracted to a much simpler, but mathematically equivalent, representation, the required analysis is performed using the abstraction, and analytic conclusions are then mapped back to the original network and interpreted there. We identify an important class of complex networks which admit vulnerability-preserving, finite state abstractions, provide efficient algorithms for computing these abstractions, and offer a vulnerability analysis methodology which combines finite state network representations with formal analytics to enable rigorous vulnerability analysis for networks of real-world scale and complexity. The considerable potential of the method is demonstrated through a case study involving a realistic electric power grid model.

We now demonstrate that the proposed approach to analyzing complex network dynamics can also be applied to biological and social systems. Consider first a biological example. Many aspects of the physiology of living organisms oscillate with a period of approximately 24 hours, corresponding to the duration of a day, and the molecular basis for this circadian rhythm has been quantified in several organisms. For instance, a useful model for the gene regulatory network responsible for circadian rhythm in *Drosophila melanogaster* (fruit fly) is [20]:

$$\frac{dM_P}{dt} = v_{sP} \frac{K_{IP}^n}{K_{IP}^n + C_N^n} - v_{mP} \frac{M_P}{K_{mP} + M_P} - k_d M_P,$$

$$\frac{dP_0}{dt} = k_{sP} M_P - V_{1P} \frac{P_0}{K_{1P} + P_0} + V_{2P} \frac{P_1}{K_{2P} + P_1} - k_d P_0,$$

$$\frac{dP_1}{dt} = V_{1P} \frac{P_0}{K_{1P} + P_0} - V_{2P} \frac{P_1}{K_{2P} + P_1} - V_{3P} \frac{P_1}{K_{3P} + P_1}$$
$$+ V_{4P} \frac{P_2}{K_{4P} + P_2} - k_d P_1,$$

$$\frac{dP_2}{dt} = V_{3P} \frac{P_1}{K_{3P} + P_1} - V_{4P} \frac{P_2}{K_{4P} + P_2} - k_3 P_2^2 + k_4 C$$
$$- v_{dP} \frac{P_2}{K_{dP} + P_2} - k_d P_2,$$

$$\frac{dC}{dt} = k_3 P_2^2 - k_4 C - k_1 C + k_2 C_N - k_{dC} C,$$

$$\frac{dC_N}{dt} = k_1 C - k_2 C_N - k_{dC} C,$$

where $M_P$, $P_0$, $P_1$, $P_2$, C, and $C_N$ are state variables corresponding to the concentrations of the constituents of the circadian rhythm gene network, $v_{sP}$ is an exogenous (control) input signal associated with the light-dark cycle of the environment, and all other terms are constant model parameters.

As is evident from Definition 2.7, a differentially flat system possesses (flat) outputs, equal in number to the number of inputs, which permit the system states and inputs to be recovered through algebraic manipulation of these outputs and their time derivatives. In the case of *Drosophila* circadian rhythm, $C_N$ is the flat output. To see this, note that C and its time derivatives can be obtained from the sixth equation through manipulation of $C_N$ and its derivatives. These terms, in turn, permit $P_2$ (and its derivatives) to be obtained from the fifth equation, and continuing in this way up the "chain" of equations gives all of the states and the input $v_{sP}$. Thus the system states and input can be obtained from knowledge of $C_N$ and its derivatives, proving that the above gene network model for *Drosophila* circadian rhythm is differentially flat. This, in turn, implies that the model admits a finite state bisimulation (Theorem 3). We have applied Algorithm 4.1 to this finite state model and identified gene network vulnerabilities which are consistent with those discussed in [21].

Consider next the phenomenon of social movements, that is, large, informal groupings of individuals and/or organizations focused on a particular issue, for instance of political, social, economic, or religious significance [e.g., 22]. Given the importance of social movements and the desire to understand their emergence and growth, numerous mathematical representations have been proposed to characterize their dynamics. For example, [23] suggests a model in which each individual in a population of interest can be in one of three states – member (of the movement), potential member, and ex-member – and interactions between individuals can lead to transitions between these affiliation states (e.g., potential members can be "persuaded" to become members). In particular, [23] proposes the following model for social movement dynamics:

$$dP/dt = \Lambda - \beta PM + \delta_1 E,$$
$$dM/dt = \beta PM - \delta_2 ME - \delta_3 M,$$
$$dE/dt = \delta_2 ME + \delta_3 M - \delta_1 E,$$

where P, M, and E denote the fractions of potential members, members, and ex-members in the population, $\Lambda$ can be interpreted to be the system's input, and $\beta$, $\delta_1$, $\delta_2$, $\delta_3$ are nonnegative constants related to the probabilities of individuals undergoing the various state transitions. It is worth noting that this model is shown in [23] to provide a good description for the growth of real world social movements.

This model for social dynamics is differentially flat with flat output E. To see this, observe that M and its time derivatives can be obtained from the third equation through manipulation of E and its derivatives. These terms, in turn, permit P (and its derivatives) to be obtained from the second equation. Finally, knowledge of P, M, E, and their derivatives allows the input $\Lambda$ to be recovered from the first equation. Thus all of the system states as well as the input can be obtained from knowledge of E and its derivatives, which shows that the social movement model is differentially flat. This, in turn, implies that the model admits a finite state bisimilar abstraction (see Theorem 3).